\documentclass[12pt]{article}
\usepackage{amssymb,amsmath,amscd, amsthm,stmaryrd,verbatim, mathrsfs}
\numberwithin{equation}{section}

\begin{document}
\baselineskip=14pt

\newcommand{\la}{\langle}
\newcommand{\ra}{\rangle}
\newcommand{\psp}{\vspace{0.4cm}}
\newcommand{\pse}{\vspace{0.2cm}}
\newcommand{\ptl}{\partial}
\newcommand{\dlt}{\delta}
\newcommand{\sgm}{\sigma}
\newcommand{\al}{\alpha}
\newcommand{\be}{\beta}
\newcommand{\G}{\Gamma}
\newcommand{\gm}{\gamma}
\newcommand{\vs}{\varsigma}
\newcommand{\Lmd}{\Lambda}
\newcommand{\lmd}{\lambda}
\newcommand{\td}{\tilde}
\newcommand{\vf}{\varphi}
\newcommand{\yt}{Y^{\nu}}
\newcommand{\wt}{\mbox{wt}\:}
\newcommand{\rd}{\mbox{Res}}
\newcommand{\ad}{\mbox{ad}}
\newcommand{\stl}{\stackrel}
\newcommand{\ol}{\overline}
\newcommand{\ul}{\underline}
\newcommand{\es}{\epsilon}
\newcommand{\dmd}{\diamond}
\newcommand{\clt}{\clubsuit}
\newcommand{\vt}{\vartheta}
\newcommand{\ves}{\varepsilon}
\newcommand{\dg}{\dagger}
\newcommand{\tr}{\mbox{Tr}}
\newcommand{\ga}{{\cal G}({\cal A})}
\newcommand{\hga}{\hat{\cal G}({\cal A})}
\newcommand{\Edo}{\mbox{End}\:}
\newcommand{\for}{\mbox{for}}
\newcommand{\kn}{\mbox{ker}}
\newcommand{\Dlt}{\Delta}
\newcommand{\rad}{\mbox{Rad}}
\newcommand{\rta}{\rightarrow}
\newcommand{\mbb}{\mathbb}
\newcommand{\lra}{\Longrightarrow}
\newcommand{\X}{{\cal X}}
\newcommand{\Y}{{\cal Y}}
\newcommand{\Z}{{\cal Z}}
\newcommand{\U}{{\cal U}}
\newcommand{\V}{{\cal V}}
\newcommand{\W}{{\cal W}}
\newcommand{\sta}{\theta}
\setlength{\unitlength}{3pt}
\newcommand{\msr}{\mathscr}
\newcommand{\wht}{\widehat}
\newcommand{\mfk}{\mathfrak}

\begin{center}{\large \bf Congruence Classes of Supporting  the Erd\"{o}s-Straus \\\pse Conjecture I: Tame Solutions} \footnote {2010 Mathematical Subject
Classification. Primary 11D68: Secondary 11D85, 11A67, 11B75}
\end{center}

\vspace{0.2cm}

\begin{center}{\large  Xiaoping Xu}\end{center}

\pse

\begin{center}{
 HLM, Institute of Mathematics, Academy of Mathematics \& System
Sciences\\ Chinese Academy of Sciences, Beijing 100190, P.R. China
\\ \& School of Mathematics, University of Chinese Academy of Sciences,\\ Beijing 100049, P.R. China}\end{center}

\begin {abstract}
\quad

In 1948, Erd\"{o}s and Straus formulated a conjecture : for any positive integer $n>2$, there exist positive integers $n_1,n_2$ and $n_3$ such that
\begin{equation}\frac{4}{n}=\frac{1}{n_1}+\frac{1}{n_2}+\frac{1}{n_3},\nonumber\end{equation}
which is still open. It is known that one only needs to prove the conjecture for any prime number $n$ such that $n\equiv 1\;(\mbox{mod}\;24)$.
If $n=24m+1$ and $n_1\leq n_2,n_3$, then $n_1=6m+k$ with $1\leq k\leq 12m$. A solution $(n_1,n_2,n_3)$ of the above equation is called a {\it tame solution} if $n_2$ and $n_3$ are factors of $(6m+k)(24m+1)$. We call $n=24m+1$ {\it wild} if it does not have any tame solution. In this paper, we
derive twenty-eight families of tame solutions of the above equation. Numeric test shows that they cover all the tame primes.

 \vspace{0.3cm}
\noindent{\it Keywords}:\hspace{0.3cm} Erd\"{o}s-Straus conjecture; Egyptian fraction; congruence class; tame solution; wild solution; wild prime.
\end{abstract}

\section {Introduction}

Ancient Egyptians used sums of unit fractions (whose numerators are 1) to express fractions due to their ways of distributing food.
For instance,  $5/8$ was interpreted by ancient Egyptians as distributing five pancakes  fairly among eight people. They cut first four
pancakes into halves and then cut the last one into eight equal pieces. So each person got the same share: $1/2+1/8$ pancakes. This amazingly interpreted the mathematical equation
\begin{equation}\frac{5}{8}=\frac{1}{2}+\frac{1}{8}.\end{equation}
So an Egyptian fraction is a sum of distinct unit fractions. By repeatedly applying the simple fact
\begin{equation}\frac{1}{k}=\frac{1}{k+1}+\frac{1}{k(k+1)},\end{equation}
one can easily prove that any fraction is a Egyptian fraction. However,  it is very difficult to determine if a fraction can be expressed as a sum of fixed number of unit fractions (cf. \cite{ET} for an excellent exposition and extensive references). For example, it is difficult to know if a fraction can be written as a sum of two unit fractions (e.g., cf. \cite{HV1,HV2,HV3,Jc2, Ps}). In 1948, Erd\"{o}s and Straus formulated a conjecture : for any positive integer $n>2$, there exist positive integers $n_1,n_2$ and $n_3$ such that
\begin{equation}\frac{4}{n}=\frac{1}{n_1}+\frac{1}{n_2}+\frac{1}{n_3},\end{equation}
which is still open up to now. For each $n$, the number of solution can go to infinity as $n$ does. Elscholtz and Tao \cite{ET} found excellent bounds for it.

Mordell \cite{Ml} proved that the conjecture holds for positive integers in  834 congruence classes modulo 840. Terzi \cite{Tdj} used computer to verify that the conjecture holds for positive integers in  the congruence classes modulo 120120 except 198 classes. Kotsireas \cite{Ki}  verified the conjecture for every $n<10^{10}$. There are other interesting partial results on the conjecture or related works (e.g., cf. \cite{Bk, BL, Ht, IW, Ld, Ms, Ps, Rl, Sjw1, Sjw2, Sjw3, Sj, ST, Tdj, Vr, Ww}).

It is known that one only needs to prove the conjecture for any prime number $n$ such that $n\equiv 1\;(\mbox{mod}\;24)$, which will be reviewed in next section. Suppose that $n=24m+1$ and $n_1\leq n_2,n_3$. Then $n_1=6m+k$ with $1\leq k\leq 12m$ and
\begin{equation}\frac{4}{24m+1}=\frac{1}{6m+k}+\frac{4k-1}{(6m+k)(24m+1)}.\end{equation}
If
\begin{equation}\frac{1}{n_2}=\frac{\Im_1}{(6m+k)(24m+1)}\quad\mbox{and}\quad \frac{1}{n_2}=\frac{\Im_2}{(6m+k)(24m+1)}\end{equation}for some positive integers $\Im_1$ and $\Im_2$ such that $\Im_1+\Im_2=4k-1$, we call the triple $(n_1,n_2,n_3)$ a {\it tame solution} of the Erd\"{o}s-Straus equation (1.3) for the positive integer $n=24m+1$. Moreover, we call $n=24m+1$ {\it wild} if it does not have any tame solution. Based on our definition in the earlier version of this work \cite{Xx1}, Howerton \cite{Hb}  showed that there are only fourteen wild primes of the form $24m+1\leq 2.4\times 10^{11}$. In this paper, we
derive twenty-eight families of tame solutions of Erd\"{o}s-Straus equation (1.3). Numeric test shows that they cover all the tame primes.
Furthermore, we call $\Im_1$ and $\Im_2$ the {\it numerator summands} for the tame solution $(n_1,n_2,n_3)$. The paper is organized as follows.\pse

In Section 2, we present some basic known facts about the Erd\"{o}s-Straus equation and prove that a tame solution $(n_1,n_2,n_3)$ of the equation for a prime $n=24m+1$ must satisfy (1.4), (1.5) and $\Im_1,\Im_2$ are factors of $6m+k$.
In Section 3, we solve the equation for a prime of the form $24m+1$ with $\Im_1=2(2j+1)$ and $\Im_2=4\ell+1$.  In Section 4, we solve the equation for a prime of the form $24m+1$ with $\Im_1=4(2j+1)$ and $\Im_2=4\ell+3$. In Section 5, we solve the equation for a prime of the form $24m+1$ with $\Im_1=2^{r+3}(2j+1)$ and $\Im_2=4\ell+3$.

\section{Basics of the  Erd\"{o}s-Straus Equation}

In this section, we list some basic known facts about the Erd\"{o}s-Straus equation (1.3). Then we prove that the numerator summands of the tame solution of the equation for a prime $n=24m+1$ must be factors of  $6m+k$.

As shown in \cite{Sj}, a simple calculation
\begin{equation}\frac{4}{2k}=\frac{1}{k}+\frac{1}{k+1}+\frac{1}{k(k+1)},\end{equation}
\begin{equation}\frac{4}{3k}=\frac{1}{3k}+\frac{1}{k+1}+\frac{1}{k(k+1)},\end{equation}
and
\begin{equation}\frac{4}{3k-1}=\frac{1}{3k-1}+\frac{1}{k}+\frac{1}{k(3k-1)}\end{equation}
shows that Erd\"{o}s-Straus equation holds for all $n$ except those $n\equiv 1\; (\mbox{mod}\;6)$. Throughout this paper, we denote by $\mbb N$ the set of nonnegative integers. So we only need to solve (1.3) for $n\in 1+6\mbb N$. Further calculations
\begin{equation}\frac{4}{4k-1}=\frac{1}{k}+\frac{1}{k(4k-1)+1}+\frac{1}{k(4k-1)(k(4k-1)+1)}\end{equation}
and
\begin{equation}\frac{4}{24k-11}=\frac{1}{6k-2}+\frac{1}{(3k-1)(24k-11)}+\frac{1}{(6k-2)(24k-11)}\end{equation}
show that Erd\"{o}s-Straus equation holds for all $n$ except those $n\equiv 1\; (\mbox{mod}\;24)$. This fact is known to some people (e.g., cf. \cite{Ms}). If $(n_1,n_2,n_3)$ is a solution of the  Erd\"{o}s-Straus equation (1.3) for $n$ and $k$ is another positive integer, then
\begin{equation}\frac{4}{nk}=\frac{1}{n_1k}+\frac{1}{n_2k}+\frac{1}{n_3k}\end{equation}
naturally holds. So it is enough to solve the equation for any prime number $n$.  From now on, we always assume that \begin{equation}n=24m+1\quad\mbox{is a prime and}\;m\in\mbb N.\end{equation}
Moreover, we can assume a tame solution $(n_1,n_2,n_3)$ of (1.3) for $n$ satisfying $n_1\leq n_2,n_3$. Then we have
\begin{equation}\frac{4}{3(24m+1)}\leq\frac{1}{n_1}\leq \frac{4}{24m+1}.
\end{equation}
Equivalently
\begin{equation}6m+1\leq n_1\leq 18m.\end{equation}
Thus
\begin{equation}n_1=6m+k\;\;\mbox{with}\;\; 1\leq k\leq 12m.\end{equation}
In particular,
\begin{equation}\frac{4}{24m+1}=\frac{1}{6m+k}+\frac{4k-1}{(6m+k)(24m+1)}\end{equation}
and
\begin{equation}\frac{1}{n_2}=\frac{\Im_1}{(6m+k)(24m+1)}\quad\mbox{and}\quad \frac{1}{n_2}=\frac{\Im_2}{(6m+k)(24m+1)}\end{equation}for some positive integers $\Im_1$ and $\Im_2$ such that
\begin{equation}\Im_1+\Im_2=4k-1. \end{equation}

Denote
\begin{equation}\Im=\mbox{l.c.m}(\Im_1,\Im_2),\end{equation}
the least common multiple of $\Im_1$ and $\Im_2$. Let $c\in \mbb Z$ such that $m+c\Im>0$. Observe
\begin{equation}(6(m+c\Im)+k)(24(m+c\Im)+1)=(6m+k)(24m+1)+6c\Im(48m+24c\Im+4k+1). \end{equation}
By Assumption (2.12),
\begin{equation}(6m+k)(24m+1)=n_2\Im_1=n_3\Im_2.\end{equation}
The above two expressions imply that \begin{equation}\frac{\Im_1}{(6(m+c\Im)+k)(24(m+c\Im)+1)}\quad\mbox{and}\quad \frac{\Im_2}{(6(m+c\Im)+k)(24(m+c\Im)+1)}\end{equation}are unit fractions. Moreover,
\begin{equation}4(6(m+c\Im)+k)-(24(m+c\Im)+1)=4(6m+1)-(2m+1)=4k-1=\Im_1+\Im_2.\end{equation}
Therefore,
\begin{eqnarray}\frac{4}{24(m+c\Im)+1}&=&\frac{1}{6(m+c\Im)+k}+\frac{\Im_1}{(6(m+c\Im)+k)(24(m+c\Im)+1)}\nonumber\\&&+\frac{\Im_1}{(6(m+c\Im)+k)(24(m+c\Im)+1)}
\end{eqnarray}
is a solution of the Erd\"{o}s-Straus equation (1.3); that is,
\begin{equation}\{24(m+c\Im)+1\mid c\in\mbb Z\;\mbox{such that}\;m+c\Im>0\}\quad\mbox{are tame numbers}.\end{equation}
This partially explains why Mordell \cite{Ml} and Terzi \cite{Tdj} had their modulo conditions. \psp

Note that
\begin{equation}4k-1\leq 48m-1\end{equation}
by (2.10). Suppose $\Im_i\not|(6m+k)$. Since $24m+1$ is a prime, we must have $(24m+1)|\Im_i$. So the above equation yields
\begin{equation}\Im_i=24m+1.\end{equation}
Without loss of generality, we may assume $i=2$. Then
\begin{equation}\Im_1=4k-1-\Im_2=4k-24m-2\leq 24m-2.\end{equation}
Thus
\begin{equation}\Im_1|(6m+k).\end{equation}
According to (2.23), we can write
\begin{equation}\Im_1=4j+2\qquad\mbox{with}\;\;j\in\mbb N.\end{equation}
Based on (2.11) and (2.13), we have
\begin{equation}\Im_1+24m+1=\Im_1+\Im_2=4k-1=4(6m+k)-(24m+1).\end{equation}
By (2.24),
\begin{equation}\Im_1|[2(24m+1)]\lra (2j+1)|(24m+1).\end{equation}
If $j=0$, $24m+3=4k-1\lra k=6m+1$ is odd, which contradicts (2.24) and (2.25).
So $j>0$. Then (2.23) and (2.25) show that $2j+1\geq 3$ is a proper factor of $24m+1$. This contradicts the assumption that $24m+1$ is a prime.
\psp

{\bf Theorem 2.1}\quad {\it For a prime $n$ of the form $24m+1$, its tame solution $(n_1,n_2,n_3)$ must satisfies (2.10)-(2.13) and}
\begin{equation}\Im_1|(6m+k),\quad \Im_2|(6m+k).\end{equation}\pse

Since
\begin{equation}\Im_1+\Im_2=4k-1=4(6m+k)-(24m+1)\end{equation}
and $n=24m+1$ is a prime, we must have
\begin{equation}\mbox{g.c.d}(\Im_1,\Im_2)=1.\end{equation}

In the rest of this paper, we always derive a necessary condition of $m$ on the tameness of a prime $24m+1$, but the condition is also sufficient for any $m$ satisfying the condition such that $24m+1$ is tame.

\section{Cases with the $\Im_1=2(2j+1)$ and $\Im_2=4\ell+1$}

 First we have
 \begin{equation}k=j+\ell+1\end{equation}by (2.13).
Moreover,
\begin{equation}6m+k=6m+j+\ell+1.\end{equation}
According to the firs expression in (2.28),
\begin{equation}2|(6m+j+\ell+1).\end{equation}So
\begin{equation}j+\ell\equiv 1\quad(\mbox{mod}\;2).\end{equation}

\subsection{Case  $j=2s$ and $\ell=6t+1$}

In this case,
 \begin{equation}\Im_1=2(4s+1),\quad \Im_2=24t+5\end{equation} and
 \begin{equation}6m+j+\ell+1=2(3(m+t)+s+1))\end{equation}
According to the firs expression in (2.28),
\begin{equation}(4s+1)|(3(m+t)+s+1);\end{equation}
that is,
\begin{equation}(3(m+t)+s+1)\equiv 0\quad(\mbox{mod}\;4s+1).\end{equation}
So \begin{equation}3(m+t)+s+1\equiv 4s+1\quad(\mbox{mod}\;12s+1).\end{equation}
Equivalently,
 \begin{equation}3(m+t)\equiv 3s\quad(\mbox{mod}\;4s+1).\end{equation}\pse

 {\it Subcase 1}. $s\not\equiv 2\;\;(\mbox{mod}\;3)$.\psp

 In this subcase, $4s+1\not \equiv 2\;\;(\mbox{mod}\;3)$ and
 the above equation implies
 \begin{equation}m+t\equiv s\quad(\mbox{mod}\;4s+1);\end{equation}
 that is,
 \begin{equation}m+t=s+c(4s+1)\qquad\mbox{for some}\;\;c\in\mbb N.\end{equation}
 So
 \begin{equation}6m+j+\ell+1=2(3c+1)(4s+1).\end{equation}
According to the second expression in (2.28) and (2.30),
 \begin{equation}(24t+5)|(3c+1)\lra c=16t+3+d(24t+5)\quad \for\;d\in\mbb N.\end{equation}Thus
 \begin{equation}m=s-t+(4s+1)[16t+3+d(24t+5)]\qquad\mbox{for some}\;\;c\in\mbb N.\end{equation}
 and
  \begin{equation}6m+j+\ell+1=2(3d+2)(4s+1)(24t+5).\end{equation}
Expression (2.11) becomes
\begin{eqnarray}&&\frac{4}{24m+1}=\frac{1}{2(3d+2)(4s+1)(24t+5)}\nonumber\\ & &+\frac{2(4s+1)+(24t+5)}{2(3d+2)(4s+1)(24t+5)(24m+1)}\nonumber\\&=&
\frac{1}{2(3d+2)(4s+1)(24t+5)}+\frac{1}{(3d+2)(24t+5)(24m+1)}\nonumber\\ & &+\frac{1}{2(3d+2)(4s+1)(24m+1)}.\end{eqnarray}\pse

 {\it Subcase 2}. $s=3s_1+2$ with $s_1\in\mbb N$.\psp

In this subcase,
\begin{equation}4s+1=12s_1+9=3(4s_1+3).\end{equation}So (3.10) is equivalent to
\begin{equation}3(m+t)\equiv 9s_1+6\quad(\mbox{mod}\;3(4s_1+3),\end{equation}\pse
 which yields
\begin{equation}m+t\equiv 3s_1+2\quad(\mbox{mod}\;4s_1+3);\end{equation}
 that is,
 \begin{equation}m+t=3s_1+2+c(4s_1+3)\qquad\mbox{for some}\;\;c\in\mbb N.\end{equation}
 So
 \begin{equation}6m+j+\ell+1=6(c+1)(4s+1).\end{equation}
According to the second expression in (2.28) and (2.30),
 \begin{equation}(24t+5)|(c+1)\lra c=-1+(d+1)(24t+5)\quad \for\;d\in\mbb N.\end{equation}Thus
 \begin{equation}m=-s_1-t-1+(d+1)(4s_1+3)(24t+5)\end{equation}
 \begin{equation}\Im_1=6(4s_1+3)\end{equation}
  and
  \begin{equation}6m+j+\ell+1=6(d+1)(4s_1+3)(24t+5).\end{equation}
Expression (2.11) becomes
\begin{eqnarray}&&\frac{4}{24m+1}=\frac{1}{6(d+1)(4s_1+3)(24t+5)}+\frac{6(4s_1+3)+(24t+5)}{6(d+1)(4s_1+3)(24t+5)(24m+1)}\nonumber\\&=&
\frac{1}{6(d+1)(4s_1+3)(24t+5)}+\frac{1}{(d+1)(24t+5)(24m+1)}\nonumber\\ & &+\frac{1}{6(d+1)(4s_1+3)(24m+1)}.\end{eqnarray}\pse

{\bf Theorem 3.1}\quad {\it For any positive integer $m$ in (3.15), we have the tame solution (3.17) of the Erd\"{o}s-Straus equation.
 When $m$ is in (3.24), we have the tame solution (3.27) of the Erd\"{o}s-Straus equation.
 }\psp

\subsection{Case  $j=6s$ and $\ell=6t+3$}

In this case,
 \begin{equation}\Im_1=2(12s+1),\quad \Im_2=24t+13\end{equation}by (3.1), and
 \begin{equation}6m+j+\ell+1=2(3(m+s+t)+2))\end{equation}
According to the firs expression in (2.28),
\begin{equation}(12s+1)|(3(m+s+t)+2);\end{equation}
that is,
\begin{equation}3(m+s+t)+2\equiv 0\quad(\mbox{mod}\;12s+1).\end{equation}
So \begin{equation}3(m+s+t)+2\equiv 2(12s+1)\quad(\mbox{mod}\;12s+1).\end{equation}
Equivalently,
 \begin{equation}3(m+t)\equiv 21s\quad(\mbox{mod}\;12s+1).\end{equation}Hence
 \begin{equation}m+t\equiv 7s\quad(\mbox{mod}\;12s+1).\end{equation}
Under this condition,
 \begin{equation}m+t=7s+c(12s+1)\qquad\mbox{for some}\;c\in\mbb N\end{equation}
 and
 \begin{equation}6m+j+\ell+1=2(3c+2)(12s+1).\end{equation}
According to the second expression in (2.28)  and (2.30),
 \begin{equation}(24t+13)|(3c+2)\lra c=16t+8+d(24t+13)\quad \for\;d\in\mbb N.\end{equation}Now
  \begin{eqnarray}m&=&7s-t+(16t+8+d(24t+13))(12s+1)\end{eqnarray} and
 \begin{equation}6m+j+\ell+1=2(3d+2)(12s+1)(24t+13).\end{equation}
  Expression (2.11) becomes
\begin{eqnarray}&&\frac{4}{24m+1}=\frac{1}{2(3d+2)(12s+1)(24t+13)}+\frac{2(12s+1)+(24t+13)}{2(3d+2)(12s+1)(24t+13)(24m+1)}\nonumber\\&=&
\frac{1}{2(3d+2)(12s+1)(24t+13)}+\frac{1}{(3d+2)(24t+13)(24m+1)}\nonumber\\ & &+\frac{1}{2(3d+2)(12s+1)(24m+1)}.\end{eqnarray}\pse

{\bf Theorem 3.2}\quad {\it For any positive integer $m$ in (3.38), we have the tame solution (3.40) of the Erd\"{o}s-Straus equation.
 }\psp

\subsection{Case  $j=6s$ and $\ell=6t+5$}

In this case,
 \begin{equation}\Im_1=2(12s+1),\quad \Im_2=24t+21=3(8t+7)\end{equation}by (4.1), and
 \begin{equation}6m+j+\ell+1=6(m+s+t+1)\end{equation}
According to the firs expression in (2.28),
\begin{equation}(12s+1)|(m+s+t+1);\end{equation}
that is,
\begin{equation}m= -s-t-1+c(12s+1)\quad\mbox{for some}\;\;c\in\mbb N.\end{equation}
Moreover,
 \begin{equation}6m+j+\ell+1=6c(12s+1).\end{equation}
According to the second expression in (2.28) and (2.30),
 \begin{equation}(24t+21)|[6c(12s+1)]\lra (8t+7)|c\lra c=(d+1)(8t+7)\end{equation}
for some $d\in\mbb N$. Thus
 \begin{equation}m= -s-t-1+(d+1)(12s+1)(8t+7)\end{equation}and
 \begin{equation}6m+j+\ell+1=6(d+1)(12s+1)(8t+7).\end{equation}
 Expression (2.11) becomes
\begin{eqnarray}&&\frac{4}{24m+1}=\frac{1}{6(d+1)(12s+1)(8t+7)}\nonumber\\ & &+\frac{2(12s+1)+3(8t+7)}{6(d+1)(12s+1)(8t+7)(24m+1)}\nonumber\\&=&
\frac{1}{6(d+1)(12s+1)(8t+7)}+\frac{1}{3(d+1)(8t+7)(24m+1)}\nonumber\\ & &+\frac{1}{2(d+1)(12s+1)(24m+1)}.\end{eqnarray}\pse

{\bf Theorem 3.3}\quad {\it For any positive integer $m$ in (3.47), we have the tame solution (3.49) of the Erd\"{o}s-Straus equation.}\psp

\subsection{Case  $j=6s+1$ and $\ell=6t+4$}

In this case,
 \begin{equation}\Im_1=6(4s+1),\quad \Im_2=24t+17\end{equation}by (3.1), and
 \begin{equation}6m+j+\ell+1=6(m+s+t+1)\end{equation}
According to the first expression in (2.28),
\begin{equation}(4s+1)|(m+s+t+1);\end{equation}
that is,
\begin{equation}m= -s-t-1+c(4s+1)\quad\mbox{for some}\;\;c\in\mbb N.\end{equation}
Moreover,
 \begin{equation}6m+j+\ell+1=6c(4s+1).\end{equation}
According to the second expression in (2.28) and (2.30),
 \begin{equation}(24t+17)|[6c(4s+1)]\lra (24t+17)|c\lra c=(d+1)(24t+17)\end{equation}
for some $d\in\mbb N$. Thus
 \begin{equation}m= -s-t-1+(d+1)(4s+1)(24t+17)\end{equation}and
 \begin{equation}6m+j+\ell+1=6(d+1)(4s+1)(24t+17).\end{equation}
 Expression (2.11) becomes
\begin{eqnarray}&&\frac{4}{24m+1}=\frac{1}{6(d+1)(4s+1)(24t+17)}+\frac{6(4s+1)+(24t+17)}{6(d+1)(4s+1)(24t+17)(24m+1)}\nonumber\\&=&
\frac{1}{6(d+1)(4s+1)(24t+17)}+\frac{1}{(d+1)(24t+17)(24m+1)}\nonumber\\ & &+\frac{1}{6(d+1)(4s+1)(24m+1)}.\end{eqnarray}\pse

{\bf Theorem 3.4}\quad {\it For any positive integer $m$ in (3.56), we have the tame solution (3.58) of the Erd\"{o}s-Straus equation.
 }\psp

\subsection{Case  $j=6s+2$ and $\ell=6t+3$}

In this case,
 \begin{equation}\Im_1=2(12s+5),\quad \Im_2=24t+13\end{equation}by (3.1), and
 \begin{equation}6m+j+\ell+1=6(m+s+t+1)\end{equation}
According to the firs expression in (2.28),
\begin{equation}(12s+5)|(m+s+t+1);\end{equation}
that is,
\begin{equation}m= -s-t-1+c(12s+5)\quad\mbox{for some}\;\;c\in\mbb N.\end{equation}
Moreover,
 \begin{equation}6m+j+\ell+1=6c(12s+5).\end{equation}
According to the second expression in (2.28) and (2.30),
 \begin{equation}(24t+13)|c\lra c=(d+1)(24t+13)\end{equation}
for some $d\in\mbb N$. Thus
 \begin{equation}m= -s-t-1+(d+1)(12s+5)(24t+13)\end{equation}and
 \begin{equation}6m+j+\ell+1=6(d+1)(12s+5)(24t+13).\end{equation}
 Expression (2.11) becomes
\begin{eqnarray}&&\frac{4}{24m+1}=\frac{1}{6(d+1)(12s+5)(24t+13)}\nonumber\\ & &+\frac{2(12s+5)+(24t+13)}{6(d+1)(12s+5)(24t+13)(24m+1)}\nonumber\\&=&
\frac{1}{6(d+1)(12s+5)(24t+13)}+\frac{1}{3(d+1)(24t+13)(24m+1)}\nonumber\\ & &+\frac{1}{6(d+1)(12s+5)(24m+1)}.\end{eqnarray}\pse

{\bf Theorem 3.5}\quad {\it For any positive integer $m$ in (3.65), we have the tame solution (3.67) of the Erd\"{o}s-Straus equation.
 }\psp

\subsection{Case  $j=6s+3$ and $\ell=6t$}

In this case,
 \begin{equation}\Im_1=2(12s+7),\quad \Im_2=24t+1\end{equation}by (3.1), and
 \begin{equation}6m+j+\ell+1=2(3(m+s+t)+2)).\end{equation}
According to the firs expression in (2.28),
\begin{equation}(12s+7)|(3(m+s+t)+2);\end{equation}
that is,
\begin{equation}3(m+s+t)+2\equiv 0\quad(\mbox{mod}\;12s+7).\end{equation}
So \begin{equation}3(m+s+t)+2\equiv 2(12s+7)\quad(\mbox{mod}\;12s+7).\end{equation}
Equivalently,
 \begin{equation}3(m+t)\equiv 21s+12\quad(\mbox{mod}\;12s+7).\end{equation}So
 \begin{equation}m+t\equiv 7s+4\quad(\mbox{mod}\;12s+7).\end{equation}
Under this condition,
 \begin{equation}m+t=7s+4+c(12s+7)\qquad\mbox{for some}\;c\in\mbb N\end{equation}
 and
 \begin{equation}6m+j+\ell+1=2(3c+2)(12s+7)).\end{equation}
The second expression in (2.28) and (2.30) imply
 \begin{equation}(24t+1)|(3c+2)\lra c=16t+d(24t+1)\quad \for\;d\in\mbb N.\end{equation}Now
  \begin{equation}m=7s-t+4+(16t+d(24t+1))(12s+7).\end{equation}and
\begin{equation}6m+j+\ell+1=2(3d+2)(12s+7)(24t+1).\end{equation}
Expression (2.11) becomes
\begin{eqnarray}&&\frac{4}{24m+1}=\frac{1}{2(3d+2)(12s+7)(24t+1)}\nonumber\\ & &+\frac{2(12s+7)+(24t+1)}{2(3d+2)(12s+7)(24t+1)(24m+1)}\nonumber\\&=&
\frac{1}{2(3d+2)(12s+7)(24t+1)}+\frac{1}{(3d+2)(24t+1)(24m+1)}\nonumber\\ & &+\frac{1}{2(3d+2)(12s+7)(24m+1)}.\end{eqnarray}\pse

{\bf Theorem 3.6}\quad {\it For any positive integer $m$ in (3.78), we have the tame solution (3.80) of the Erd\"{o}s-Straus equation.
 }\psp

\subsection{Case  $j=6s+3$ and $\ell=6t+2$}

In this case,
 \begin{equation}\Im_1=2(12s+7),\quad \Im_2=24t+9=3(8t+3)\end{equation}by (3.1), and
 \begin{equation}6m+j+\ell+1=6(m+s+t+1)\end{equation}
According to the firs expression in (2.28),
\begin{equation}(12s+7)|(m+s+t+1);\end{equation}
that is,
\begin{equation}m+s+t+1\equiv 0\quad(\mbox{mod}\;12s+7).\end{equation}
So
 \begin{equation}m+t\equiv 11s+6\quad(\mbox{mod}\;12s+7).\end{equation}
Under this condition,
 \begin{equation}m+t=11s+6+c(12s+7)\qquad\mbox{for some}\;c\in\mbb N\end{equation}
 and
 \begin{equation}6m+j+\ell+1=6(c+1)(12s+7).\end{equation}
The second expression in (2.28) and (2.30) imply
\begin{equation}(8t+3)|(c+1).\end{equation}Thus
\begin{equation}
c=(d+1)(8t+3)-1\quad \for\;\;0<d\in\mbb N.\end{equation}Now
  \begin{equation}m=-t-s-1+(d+1)(12s+7)(8t+3)\end{equation}and
   \begin{equation}6m+j+\ell+1=6(d+1)(12s+7)(8t+3).\end{equation}
Expression (2.11) becomes
\begin{eqnarray}\frac{4}{24m+1}&=&\frac{1}{6(d+1)(12s+7)(8t+3)}\nonumber\\ & &+\frac{2(12s+7)+3(8t+3)}{6(d+1)(12s+7)(8t+3)(24m+1)}\nonumber\\&=&
\frac{1}{6(d+1)(12s+7)(8t+3)}+\frac{1}{3(d+1)(8t+3)(24m+1)}\nonumber\\ & &+\frac{1}{2(d+1)(12s+7)(24m+1)}.\end{eqnarray}\pse

{\bf Theorem 3.7}\quad {\it For any positive integer $m$ in (3.90), we have the tame solution (3.92) of the Erd\"{o}s-Straus equation.
 }\psp

\subsection{Case  $j=6s+3$ and $\ell=6t+4$}

In this case,
 \begin{equation}\Im_1=2(12s+7),\quad \Im_2=24t+17\end{equation}by (3.1), and
 \begin{equation}6m+j+\ell+1=2(3(m+s+t)+4)\end{equation}
According to the firs expression in (2.28),
\begin{equation}(12s+7)|(3(m+s+t)+4);\end{equation}
that is,
\begin{equation}3(m+s+t)+4\equiv 0\quad(\mbox{mod}\;12s+7).\end{equation}
So \begin{equation}3(m+s+t)+4\equiv 12s+7\quad(\mbox{mod}\;12s+7).\end{equation}
Equivalently,
 \begin{equation}3(m+t)\equiv 9s+3\quad(\mbox{mod}\;12s+7).\end{equation}So
 \begin{equation}m+t\equiv 3s+1\quad(\mbox{mod}\;12s+7).\end{equation}
Under this condition,
 \begin{equation}m+t=3s+1+c(12s+7)\qquad\mbox{for some}\;c\in\mbb N\end{equation}
 and
 \begin{equation}6m+j+\ell+1=2(3c+1)(12s+7).\end{equation}
According to the second expression in (2.28) and (2.30),
 \begin{equation}(24t+17)|(3c+1)\lra c=16t+11+d(24t+17)\quad \for\;d\in\mbb N.\end{equation}Now
  \begin{equation}m=3s+1-t+(16t+11+d(24t+17))(12s+7)\end{equation} and
 \begin{equation}6m+j+\ell+1=2(3d+2)(12s+7)(24t+17).\end{equation}
Expression (2.11) becomes
\begin{eqnarray}&&\frac{4}{24m+1}=\frac{1}{2(3d+2)(12s+7)(24t+17)}\nonumber\\ & &+\frac{2(12s+7)+(24t+17)}{2(3d+2)(12s+7)(24t+17)(24m+1)}\nonumber\\&=&
\frac{1}{2(3d+2)(12s+7)(24t+17)}+\frac{1}{(3d+2)(24t+17)(24m+1)}\nonumber\\ & &+\frac{1}{2(3d+2)(12s+7)(24m+1)}.\end{eqnarray}\pse

{\bf Theorem 3.8}\quad {\it For any positive integer $m$ in (3.103), we have the tame solution (3.105) of the Erd\"{o}s-Straus equation.
 }\psp

\subsection{Case  $j=6s+5$ and $\ell=6t$}

In this case,
 \begin{equation} \Im_1=2(12s+11),\quad \Im_2=24t+1\end{equation}by (3.1), and
 \begin{equation}6m+j+\ell+1=6(m+s+t+1)\end{equation}
According to the firs expression in (2.28),
\begin{equation}(12s+11)|(m+s+t+1);\end{equation}
that is,
\begin{equation}m+s+t+1\equiv 0\quad(\mbox{mod}\;12s+11).\end{equation}
So
\begin{equation}m\equiv -s-t-1\quad(\mbox{mod}\;12s+11).\end{equation}
Under this condition,
 \begin{equation}m=-s-t-1+c(12s+11)\qquad\mbox{for some}\;c\in\mbb N\end{equation}
 and
 \begin{equation}6m+j+\ell+1=6c(12s+11).\end{equation}

The second expression in (2.28)  and (2.30) imply
 \begin{equation}(24t+1)|c\lra c=(d+1)(24t+1)\quad \for\;0<d\in\mbb N.\end{equation}Now
  \begin{equation}m=-s-t-1+(d+1)(24t+1)(12s+11).\end{equation} and
 \begin{equation}6m+j+\ell+1=6(d+1)(12s+11)(24t+1).\end{equation}
Expression (2.11) becomes
\begin{eqnarray}&&\frac{4}{24m+1}=\frac{1}{6(d+1)(12s+11)(24t+1)}\nonumber\\ & &++\frac{2(12s+11)+(24t+1)}{6(d+1)(12s+11)(24t+1)(24m+1)}\nonumber\\&=&
\frac{1}{(d+1)(12s+11)(24t+1)}+\frac{1}{3(d+1)(24t+1)(24m+1)}\nonumber\\ & &+\frac{1}{6(d+1)(12s+11)(24m+1)}.\end{eqnarray}
\pse

{\bf Theorem 3.9}\quad {\it For any positive integer $m$ in (3.145), we have the tame solution (3.147) of the Erd\"{o}s-Straus equation.
 }\psp

\section{Cases $\Im_1=4j$, $j$ is odd  and  $\Im_2=4\ell+3$}

This is the case when
 \begin{eqnarray}\Im_1=4j,\quad k=j+\ell+1\quad\mbox{with}\;\;j,\ell\in \mbb N\end{eqnarray} and $\ell\geq 1$ by (2.13).
Moreover,
\begin{equation}6m+k=6m+j+\ell+1.\end{equation}
According to the first expression in (2.28),
\begin{equation}4|(6m+j+\ell+1).\end{equation}So
\begin{equation}j+\ell\equiv 1\quad(\mbox{mod}\;2).\end{equation}

\subsection{Case $m=2m_1, j=4s+1$ and $\ell=12t+3$}

Now
\begin{equation}\Im_1=4(4s+1),\quad \Im_2=48t+11\end{equation} and
\begin{equation}6m'+j'+\ell'=4(3(m_1+t)+s+1).\end{equation}
According to the first expression in (2.28),
\begin{equation}(4s+1)|(3(m_1+t)+s+1);\end{equation}
that is
\begin{equation}3(m_1+t)+s+1\equiv 0\quad(\mbox{mod}\;4s+1).\end{equation}Note
\begin{equation}3(m_1+t)+s+1\equiv 4s+1\quad(\mbox{mod}\;4s+1).\end{equation}
Equivalently,
\begin{equation}3(m_1+t)\equiv 3s\quad(\mbox{mod}\;4s+1).\end{equation}\pse

{\it Subcase 1}.   $s\not\equiv 2\;(\mbox{mod}\;3)$.\psp

 In this subcase,
\begin{equation}m_1+t\equiv s\quad(\mbox{mod}\;4s+1).\end{equation}
Hence
\begin{equation}m_1+t= s+c(4s+1)\quad\mbox{for some}\;\;c\in\mbb N.\end{equation}
Observe that
\begin{equation}6m+j+\ell;=4(3c+1)(4s+1).\end{equation}
According to the second expression in (2.28) and (2.30),
\begin{equation}(48t+11)|(3c+1).\end{equation}
Therefore,
\begin{equation}c=32t+7+d(48t+11)\qquad\mbox{for some}\;\;d\in\mbb N.\end{equation} This implies
\begin{eqnarray}m&=&2m_1=2[s-t+(32t+7+d(48t+11))(4s+1)]\end{eqnarray}
and
 \begin{equation}6m+j+\ell+1=4(3d+2)(4s+1)(48t+11).\end{equation}
Expression (2.11) becomes
\begin{eqnarray}&&\frac{4}{24m+1}=\frac{1}{4(3d+2)(4s+1)(48t+11)}\nonumber\\ & &+\frac{4(4s+1)+(48t+11)}{4(3d+2)(4s+1)(48t+11)(24m+1)}\nonumber\\&=&
\frac{1}{4(3d+2)(4s+1)(48t+11)}+\frac{1}{(3d+2)(48t+11)(24m+1)}\nonumber\\ & &+\frac{1}{4(3d+2)(4s+1)(24m+1)}.\end{eqnarray}
\pse

{\it Subcase 2}. $s\equiv 2\;(\mbox{mod}\;3)$.\psp

In this subcase, $s=3s_1+2$ for some $s_1\in\mbb N$. Then (4.10) is equivalent to
\begin{equation}m_1+t\equiv 3s_1+2\quad(\mbox{mod}\;4s_1+3).\end{equation}
Hence
\begin{equation}m_1+t=3s_1+2+c(4s_1+3)\quad\mbox{for some}\;\;c\in\mbb N.\end{equation}
Moreover,
\begin{equation}3(m_1+t)+s+1=3(m_1+t+s_1+1)=3(c+1)(4s_1+3).\end{equation}
 According to the second expression in (2.28) and (2.30),
\begin{equation}(48t+11)|(c+1).\end{equation}Thus
\begin{equation}c=(d+1)(48t+11)-1\quad\mbox{for some}\;\;d\in\mbb N.\end{equation}
Now
\begin{equation}m=-2s_1-2t-2+2(d+1)(4s_1+3)(48t+11)\end{equation}and
\begin{equation}
6m+j+\ell+1=12(d+1)(4s_1+3)(48t+11).\end{equation}
Expression (2.11) becomes
\begin{eqnarray}&&\frac{4}{24m+1}=\frac{1}{12(d+1)(4s_1+3)(48t+11)}\nonumber\\ & &+\frac{12(4s_1+3)+(48t+11)}{12(d+1)(4s_1+3)(48t+11)(24m+1)}\nonumber\\&=&
\frac{1}{12(d+1)(4s_1+3)(48t+11)}+\frac{1}{(d+1)(48t+11)(24m+1)}\nonumber\\ & &+\frac{1}{12(d+1)(4s_1+3)(24m+1)}.\end{eqnarray}\pse

{\bf Theorem 4.1}\quad {\it For any positive integer $m$ in (4.16), we have the tame solution (4.18) of the Erd\"{o}s-Straus equation.
  When $m$ is of the form (4.24), we have the tame solution (4.26) of the Erd\"{o}s-Straus equation. }\psp

\subsection{Case $m=2m_1+1, j=12s+5$ and $\ell=4t$}

In this case,
\begin{equation}\Im_1=4(12s+5),\quad \Im_2=16t+3\end{equation}by (4.1),  and
\begin{equation}6m+j+\ell+1=12m_1+6+12s+5+4t+1=4(3(m_1+s+1)+t)
\end{equation}
According to the second expression in (2.28),
\begin{equation}(16t+3)|(3(m_1+s+1)+t);\end{equation}
that is
\begin{equation}3(m_1+s+1)+t\equiv 0\quad(\mbox{mod}\;16t+3).\end{equation}Note
\begin{equation}3(m_1+s+1)+t\equiv 16t+3\quad(\mbox{mod}\;16t+11).\end{equation}
Equivalently,
\begin{equation}3(m_1+s)\equiv 15t\quad(\mbox{mod}\;16t+3).\end{equation}\pse

{\it Subcase 1}. $t\not\equiv 0\;\;(\mbox{mod}\;3)$.\psp

In this subcase, the above equation is equivalent to
\begin{equation}m_1+s\equiv 5t\quad(\mbox{mod}\;16t+3).\end{equation}
Hence
\begin{equation}m_1=5t-s+c(16t+3)\quad\mbox{for some}\;\;c\in\mbb N.\end{equation}

Observe that
\begin{equation}6m+j+\ell+1=4(3c+1)(16t+3).\end{equation}
According to the first expression in (2.28) and (2.30),
\begin{equation}(12s+5)|(3c+1).\end{equation}
Therefore,
\begin{equation}c=8s+3+d(12s+5)\qquad\mbox{for some}\;\;d\in\mbb N.\end{equation} This implies
\begin{equation}m_1=5t-s+(8s+3+d(12s+5))(16t+3)\end{equation}
\begin{equation}\lra m=2[5t-s+(8s+3+d(12s+5))(16t+3)]+1\end{equation}
and
 \begin{equation}6m+j+\ell+1=4(3d+2)(12s+5)(16t+3).\end{equation}
Expression (2.11) becomes
\begin{eqnarray}&&\frac{4}{24m+1}=\frac{1}{4(3d+2)(12s+3)(16t+3)}\nonumber\\ & &+\frac{4(12s+5)+(16t+3)}{4(3d+2)(12s+5)(16t+3)(24m+1)}\nonumber\\&=&
\frac{1}{4(3d+2)(12s+5)(16t+3)}+\frac{1}{(3d+2)(16t+3)(24m+1)}\nonumber\\ & &+\frac{1}{4(3d+2)(12s+5)(24m+1)}.\end{eqnarray}
\pse

{\it Subcase 2}. $t\equiv 0\;\;(\mbox{mod}\;3)$.\psp

In this subcase, $t=3t_1$ for some $t\in\mbb N$. Equation (3.42) becomes
\begin{equation}3(m_1+s)\equiv 45t_1\quad(\mbox{mod}\;3(16t_1+1)).\end{equation}
So
\begin{equation}m_1+s\equiv 15t_1\quad(\mbox{mod}\;16t_1+1);\end{equation}
that is,
\begin{equation}m_1+s= 15t_1+c(16t_1+1)\qquad\mbox{for some}\;\;c\in\mbb N.\end{equation}

Observe that
\begin{equation}6m+j+\ell+1=12(c+1)(16t_1+1).\end{equation}
According to the first expression in (2.28) and (2.30),
\begin{equation}(12s+5)|(c+1).\end{equation}
Therefore,
\begin{equation}c=(d+1)(12s+5)-1\qquad\mbox{for some}\;\;d\in\mbb N.\end{equation} This implies
\begin{equation}m_1= -s-t_1-1+(d+1)(12s+5)(16t_1+1)\end{equation}
\begin{equation}
\lra m=2[-s-t_1+(d+1)(12s+5)(16t_1+1)]-1\end{equation}
and
 \begin{equation}6m+j+\ell+1=12(d+1)(12s+5)(16t_1+1).\end{equation}
Expression (2.11) becomes
\begin{eqnarray}&&\frac{4}{24m+1}=\frac{1}{12(d+1)(12s+5)(16t_1+1)}\nonumber\\ & &+\frac{4(12s+5)+3(16t_1+1)}{12(d+1)(12s+5)(16t_1+1)(24m+1)}\nonumber\\&=&
\frac{1}{12(d+1)(12s+5)(16t_1+9)}+\frac{1}{3(d+1)(16t_1+9)(24m+1)}\nonumber\\ & &+\frac{1}{4(d+1)(12s+5)(24m+1)}.\end{eqnarray}
\pse

{\bf Theorem 4.2}\quad {\it For any positive integer $m$ in (4.39), we have the tame solution (4.41) of the Erd\"{o}s-Straus equation.
When $m$ is in (4.49), we have the tame solution (4.51) of the Erd\"{o}s-Straus equation.  }\psp

\subsection{Case $m=2m_1, j=12s+7$ and $\ell=12t+4$}

Now
\begin{equation}\Im_1=4(12s+7),\quad \Im_2=48t+19\end{equation} by (4.1),  and
\begin{equation}6m+j+\ell+1=12(m_1+s+t+1).\end{equation}
According to the first expression in (2.28),
\begin{equation}(12s+7)|(m_1+s+t+1);\end{equation}
that is
\begin{equation}m_1+s+t+1\equiv 0\quad(\mbox{mod}\;12s+7).\end{equation}
Hence
\begin{equation}m_1=-s-t-1+c(12s+7)\quad\mbox{for some}\;\;c\in\mbb N.\end{equation}

Observe that
\begin{equation}6m+j+\ell+1=12c(12s+7).\end{equation}
According to the second expression in (2.28) and (2.30),
\begin{equation}(48t+19)|c.\end{equation}
Therefore,
\begin{equation}c=(d+1)(48t+19)\qquad\mbox{for some}\;\;d\in\mbb N.\end{equation} This implies
\begin{equation}m=2m_1=2[-s-t-1+(d+1)(12s+7)(48t+19)]\end{equation}
and
 \begin{equation}6m+j+\ell+1=12(d+1)(12s+7)(48t+19).\end{equation}
Expression (2.11) becomes
\begin{eqnarray}&&\frac{4}{24m+1}=\frac{1}{12(d+1)(12s+7)(48t+19)}\nonumber\\ & &+\frac{4(12s+7)+(48t+19)}{12(d+1)(12s+7)(48t+19)(24m+1)}\nonumber\\&=&
\frac{1}{12(d+1)(12s+7)(48t+19)}+\frac{1}{3(d+1)(48t+19)(24m+1)}\nonumber\\ & &+\frac{1}{12(d+1)(12s+7)(24m+1)}.\end{eqnarray}\pse

\pse

{\bf Theorem 4.3}\quad {\it For any positive integer $m$ in (4.60), we have the tame solution (4.62) of the Erd\"{o}s-Straus equation.
  }\psp

\subsection{Case $m=2m_1+1, j=12s+1$ and $\ell=12t+4$}

Now
\begin{equation}\Im_1=4(12s+1),\quad \Im_2=48t+19\end{equation} by (4.1),  and
\begin{equation}6m+j+\ell+1=12(m_1+s+t+1).\end{equation}
According to the first expression in (2.28),
\begin{equation}(12s+1)|(m_1+s+t+1);\end{equation}
that is
\begin{equation}m_1+s+t+1\equiv 0\quad(\mbox{mod}\;12s+1).\end{equation}
Hence
\begin{equation}m_1=-s-t-1+c(12s+1)\quad\mbox{for some}\;\;c\in\mbb N.\end{equation}

Observe that
\begin{equation}6m+j+\ell+1=12c(12s+1).\end{equation}
According to the second expression in (2.28) and (2.30),
\begin{equation}(48t+19)|c.\end{equation}
Therefore,
\begin{equation}c=(d+1)(48t+19)\qquad\mbox{for some}\;\;d\in\mbb N.\end{equation} This implies
\begin{equation}m=2m_1+1=2[-s-t+(d+1)(12s+1)(48t+19)]-1\end{equation}
and
 \begin{equation}6m+j+\ell+1=12(d+1)(12s+1)(48t+19).\end{equation}
Expression (2.11) becomes
\begin{eqnarray}&&\frac{4}{24m+1}=\frac{1}{12(d+1)(12s+1)(48t+19)}\nonumber\\ & &+\frac{4(12s+1)+(48t+19)}{12(d+1)(12s+1)(48t+19)(24m+1)}\nonumber\\&=&
\frac{1}{12(d+1)(12s+1)(48t+19)}+\frac{1}{3(d+1)(48t+19)(24m+1)}\nonumber\\ & &+\frac{1}{12(d+1)(12s+1)(24m+1)}.\end{eqnarray}\pse

\pse

{\bf Theorem 4.4}\quad {\it For any positive integer $m$ in (4.70), we have the tame solution (4.72) of the Erd\"{o}s-Straus equation.
  }\psp

\subsection{Case $m=2m_1+1, j=4s+3$ and $\ell=12t+2$}

In this case,
\begin{equation}\Im_1=4(4s+3),\quad \Im_2=48t+11\end{equation}by (4.1), and
\begin{equation}6m+j+\ell+1=4(3(m_1+t+1)+s).\end{equation}

According to the first expression in (2.28),
\begin{equation}(4s+3)|(3(m_1+t+1)+s);\end{equation}
that is
\begin{equation}3(m_1+t+1)+s\equiv 0\quad(\mbox{mod}\;4s+3).\end{equation}Note
\begin{equation}3(m_1+t+1)+s\equiv 4s+3\quad(\mbox{mod}\;4s+3).\end{equation}
Equivalently,
\begin{equation}3(m_1+t)\equiv 3s\quad(\mbox{mod}\;4s+3).\end{equation}\pse

{\it Subcase 1}. $s\equiv 0\;(\mbox{mod}\;3)$.\psp

In this subcase, $s=3s_1$ for some $s_1\in\mbb N$. Moreover, the above equation becomes
\begin{equation}3(m_1+t)\equiv 9s_1\quad(\mbox{mod}\;12s_1+3).\end{equation}
So
\begin{equation}m_1+t\equiv 3s_1\quad(\mbox{mod}\;4s_1+1).\end{equation}
Thus
\begin{equation}m_1=3s_1-t+c(4s_1+1)\quad\mbox{for some}\;\;c\in\mbb N.\end{equation}
In particular,
\begin{equation}6m+j+\ell+1=12(m_1+t+s_1+1)=12(c+1)(4s_1+1).\end{equation}
The second expression in (2.28)  and (2.30) yield
\begin{equation}(48t+11)|(c+1)\lra c=(d+1)(48t+11)-1\quad\mbox{for some}\;\;d\in\mbb N.\end{equation} Equation (4.71) shows
\begin{equation}m_1=-s_1-t-1+(d+1)(4s_1+1)(48t+11).\end{equation}
\begin{equation}\lra m=2[-s_1-t+(d+1)(4s_1+1)(48t+11)]-1.\end{equation}
Furthermore,
\begin{equation}6m+j+\ell+1=12(d+1)(4s_1+1)(48t+11).\end{equation}
Therefore, (2.11) becomes
\begin{eqnarray}&&\frac{4}{24m+1}=\frac{1}{12(d+1)(4s_1+1)(48t+11)}\nonumber\\& &+\frac{12(4s_1+1)+(48t+11)}{12(d+1)(4s_1+1)(48t+11)(24m+1)}\nonumber\\&=&
\frac{1}{12(d+1)(4s_1+1)(48t+11)}+\frac{1}{(d+1)(48t+11)(24m+1)}\nonumber\\& &+\frac{1}{12(d+1)(4s_1+1)(24m+1)}.\end{eqnarray}\pse

{\it Subcase 2}. $s\not\equiv 0\;(\mbox{mod}\;3)$.\psp

By (4.68),
\begin{equation}m_1+t\equiv s\quad(\mbox{mod}\;4s+3).\end{equation}
Thus
\begin{equation}m_1\equiv s-t\quad(\mbox{mod}\;4s+3).\end{equation}
Hence
\begin{equation}m_1=s-t+c(4s+3)\quad\mbox{for some}\;\;c\in\mbb N.\end{equation}

Observe that
\begin{equation}6m+j+\ell+1=4(3c+1)(4s+3).\end{equation}
 According to
 the second expression in (2.28) and (2.30),
\begin{equation}(48t+11)|(3c+1).\end{equation}
Therefore,
\begin{equation}c=32t+7+d(48t+11)\qquad\mbox{for some}\;\;d\in\mbb N.\end{equation} This implies
\begin{equation}m_1=s-t+(32t+7+d(48t+11))(4s+3)\end{equation}
\begin{equation}\lra m=2[s-t+(32t+7+d(48t+11))(4s+3)]+1\end{equation}
and
 \begin{equation}6m+j+\ell+1=4(3d+2)(4s+3)(48t+11).\end{equation}
Expression (2.11) becomes
\begin{eqnarray}&&\frac{4}{24m+1}=\frac{1}{4(3d+2)(4s+3)(48t+11)}\nonumber\\ & &+\frac{4(4s+3)+(48t+11)}{4(3d+2)(4s+3)(48t+11)(24m+1)}\nonumber\\&=&
\frac{1}{4(3d+2)(4s+3)(48t+11)}+\frac{1}{(3d+2)(48t+11)(24m+1)}\nonumber\\ & &+\frac{1}{4(3d+2)(4s+3)(24m+1)}.\end{eqnarray}

In summary, we have:
\psp

{\bf Theorem 4.5}\quad {\it For any positive integer $m$ in (4.86), we have the tame solution (4.88) of the Erd\"{o}s-Straus equation.
When $m$ is in (4.96), we have the tame solution (4.98) of the Erd\"{o}s-Straus equation.  }\psp

\subsection{Case $m=2m_1, j=12s+11$ and $\ell=4t$}

In this case,
\begin{equation}\Im_1=4(12s+11),\quad \Im_2=16t+3\end{equation}by (4.1), and
\begin{equation}6m+j+\ell+1=4(3(m_1+s+1)+t).\end{equation}

According to the second expression in (2.28),
\begin{equation}(16t+3)|(3(m_1+s+1)+t);\end{equation}
that is
\begin{equation}3(m_1+s+1)+t\equiv 0\quad(\mbox{mod}\;16t+3).\end{equation}Note
\begin{equation}3(m_1+s+1)+t\equiv 16t+3\quad(\mbox{mod}\;16t+3).\end{equation}
Equivalently,
\begin{equation}3(m_1+s)\equiv 15t\quad(\mbox{mod}\;16t+3).\end{equation}\pse

{\it Subcase 1}. $t\not\equiv 0\;(\mbox{mod}\;3)$.\psp

In this subcase,  the above equation is equivalent to
\begin{equation}m_1+s\equiv 5t\quad(\mbox{mod}\;16t+3).\end{equation}\pse
Thus
\begin{equation}m_1=5t-s+c(16t+3)\quad\mbox{for some}\;\;c\in\mbb N.\end{equation}
In particular,
\begin{equation}6m+j+\ell+1=4(3c+1)(16t+3).\end{equation}
The first expression in (2.28)  and (2.30) yield
\begin{equation}(12s+11)|(3c+1)\lra c=8s+7+d(12s+11)\quad\mbox{for some}\;\;d\in\mbb N.\end{equation} Equation (5.91) shows
\begin{equation}m_1=5t-s+(8s+7+d(12s+11))(16t+3)\end{equation}
\begin{equation}\lra m=2[5t-s+(8s+7+d(12s+11))(16t+3)].\end{equation}Furthermore,
\begin{equation}6m+j+\ell+1=4(3d+2)(12s+11)(16t+3).\end{equation}
Therefore, (2.11) becomes
\begin{eqnarray}&&\frac{4}{24m+1}=\frac{1}{4(3d+2)(12s+11)(16t+3)}\nonumber\\& &+\frac{4(12s+11)+(16t+3)}{4(3d+2)(12s+11)(16t+3)(24m+1)}\nonumber\\&=&
\frac{1}{4(3d+2)(12s+11)(16t+3)}+\frac{1}{4(3d+2)(16t+3))(24m+1)}\nonumber\\& &+\frac{1}{4(3d+2)(12s+11)(24m+1)}.\end{eqnarray}\pse

{\it Subcase 2}. $t\equiv 0\;(\mbox{mod}\;3)$.\psp

In this subcase,  $t=st_1$ for some $t_1\in\mbb N$. Now (4.104) becomes
\begin{equation}3(m_1+s+1)\equiv 45t_1+3\quad(\mbox{mod}\;3(16t_1+1));\end{equation}
equivalently,
\begin{equation}m_1+s\equiv 15t_1\quad(\mbox{mod}\;16t_1+1);\end{equation}
Thus
\begin{equation}m_1=15t_1-s+c(16t_1+1)\quad\mbox{for some}\;\;c\in\mbb N.\end{equation}
In particular,
\begin{equation}6m+j+\ell+1=12(c+1)(16t_1+1).\end{equation}
The first expression in (2.28)  and (2.30) yield
\begin{equation}(12s+11)|(c+1)\lra c=(d+1)(12s+11)-1\quad\mbox{for some}\;\;d\in\mbb N.\end{equation} Equation (5.91) shows
\begin{equation}m_1=-s-t-1+(d+1)(12s+11)(16t_1+1)\end{equation}
\begin{equation}\lra m=2[-s-t_1-1+(d+1)(12s+11)(16t_1+1)].\end{equation}Furthermore,
\begin{equation}6m+j+\ell+1=12(d+1)(12s+11)(16t_1+1).\end{equation}
Therefore, (2.11) becomes
\begin{eqnarray}&&\frac{4}{24m+1}=\frac{1}{12(d+1)(12s+11)(16t_1+1)}\nonumber\\& &+\frac{4(12s+11)+3(16t_1+1)}{12(d+1)(12s+11)(16t_1+1)(24m+1)}\nonumber\\&=&
\frac{1}{12(d+1)(12s+11)(16t_1+1)}+\frac{1}{3(d+1)(16t_1+1))(24m+1)}\nonumber\\& &+\frac{1}{4(d+1)(12s+11)(24m+1)}.\end{eqnarray}\pse

In summary, we have:
\psp

{\bf Theorem 4.6}\quad {\it  When $m$ is in (4.110), the Erd\"{o}s-Straus equation (4.112) holds.
If $m$ is in (4.119), the Erd\"{o}s-Straus equation (4.121) holds.}

\subsection{Section Summary}

In the above sections, $m=g(x,d,y)$, $\Im_1=\eta(x)$ and $\Im_2=h(y)$ are integral polynomials, where $x=s$ or $s_1$ and $y=2t$ or $y=2t_1$.
If we change $y$ to odd positive integer, they still give a solution of the Erd\"{o}s-Straus equation. In summary, we have the following table
of solutions:

\psp
\begin{center}{\bf \large Table 4.6}\end{center}\pse
\begin{center}\begin{tabular}{|c|c|c|}\hline $m$&$\Im_1$&$\Im_2$\\\hline $2x-y+2(16y+7+d(24y+11))(4x+1)$&$4(4x+1)$&$24y+11$
\\\hline $-2x-y-2+2(d+1)(4x+3)(24y+11)$&$12(4x+3)$&$24y+11$\\\hline $5y-2x+2(8x+3+d(12x+5))(8y+3)+1$&$4(12x+5)$&$8y+3$\\\hline
$-2x-y-1+2(d+1)(12x+5)(8y+1)$&$4(12x+5)$&$3(8y+1)$\\\hline $-2x-y-2+2(d+1)(12x+7)(24y+19)$&$4(12x+7)$&$24y+19$\\\hline
$-2x-y-1+2(d+1)(12x+1)(24y+19)$&$4(12x+1)$&$24y+19$\\\hline
$-2x-y+2(d+1)(4x+1)(24y+11)-1$&$12(4x+1)$&$24y+11$\\\hline $2x-y+2(16y+7+d(24y+11))(4x+3)+1$&$4(4x+3)$&$24y+11$\\\hline
$5y+2[-x+(8x+7+d(12s+11))(8y+3)]$&$4(12x+11)$&$8y+3$\\\hline $-2x-y-2+2(d+1)(12x+11)(8y+1)$&$4(12x+11)$&$3(8y+1)$\\\hline
\end{tabular}\end{center}\psp\psp

\section{Cases $\Im_1=2^{r+3}j'$, $j'$ is odd and  $\Im_2=4\ell+3$}\psp

In this situation, $j=2^{r+1}j'$.
\begin{equation}\Im_1=2^{r+3}j',\qquad \Im_2=4\ell+3.\end{equation}
Moreover,
\begin{equation}6m+j+\ell+1=6m+2^{r+1}j'+\ell+1.\end{equation}
According to the first expression in (2.28),
\begin{equation}6m+2^{r+1}j'+\ell+1\equiv 0\quad(\mbox{mod}\;2^{r+3}j').\end{equation}
So
\begin{equation}\ell\;\;\mbox{is odd}.\end{equation}\pse

\subsection{$j=3\cdot 2^{r+1}(2s+1)$ and $\ell=3\cdot2^{r+3}t+5$}

In this case,
\begin{equation}6m+j+\ell+1=6m+3\cdot2^{r+1}(2s+1)+3\cdot2^{r+3}t+6=6[m+2^r(2s+4t+1)+1].
\end{equation}
Now (5.3) becomes,
\begin{equation}6[m+2^r(2s+4t+1)+1]\equiv 0\quad(\mbox{mod}\;2^{r+3}3(2s+1));\end{equation}
equivalently,
\begin{equation}m+2^r(2s+4t+1)+1\equiv 0\quad(\mbox{mod}\;2^{r+2}(2s+1));\end{equation}
Thus
\begin{equation}m=2^{r+2}(2s+1)c-2^r(2s+4t+1)-1\qquad\mbox{for some}\;\;0<c\in\mbb N.\end{equation}
Note
\begin{equation}\Im_2=4\ell+3=3\cdot2^{r+5}t+23,\quad 6m+j+\ell+1=3\cdot2^{r+3}(2s+1)c.\end{equation}
According to the second expression in (2.28) and (2.30),
\begin{equation}(3\cdot2^{r+5}t+23)|(3\cdot 2^{r+3}(2s+1)c)\lra(3\cdot 2^{r+5}t+23)|c \end{equation}
\begin{equation}\lra c=(d+1)(3\cdot2^{r+5}t+23)\qquad\mbox{for some}\;\;d\in\mbb N.\end{equation}
Hence
\begin{equation}m=(d+1)2^{r+3}(2s+1)(3\cdot2^{r+5}t+23)-2^r(2s+4t+1)-1\end{equation}
and
\begin{equation} 6m+j+\ell+1=3(d+1)2^{r+3}(2s+1)(3\cdot2^{r+5}t+23).\end{equation}
Therefore, (2.11) becomes
\begin{eqnarray}&&\frac{4}{24m+1}=\frac{1}{3(d+1)2^{r+3}(2s+1)(3\cdot2^{r+5}t+23)}\nonumber\\ & &+\frac{(3\cdot2^{r+3}(2s+1)+(3\cdot2^{r+5}t+23)}
{3(d+1)2^{r+3}(2s+1)(3\cdot2^{r+5}t+23)d(24m+1)}\nonumber\\&=&
\frac{1}{3(d+1)2^{r+3}(2s+1)(3\cdot2^{r+5}t+23)}+\frac{1}{(d+1)(3\cdot2^{r+5}t+23)(24m+1)}
\nonumber\\ & &+\frac{1}{3(d+1)2^{r+3}(2s+1)d(24m+1)}.\end{eqnarray}
\pse

{\bf Theorem 5.1}\quad {\it  When $m$ is in (5.12), the Erd\"{o}s-Straus equation (5.14) holds.}\psp

\subsection{$j=2^{r+1}(2s+1),\;3\not|(2s+1)$ and $\ell=3\cdot2^{r+3}t+5$}

In this case,
\begin{equation}6m+j+\ell+1=6m+2^{r+1}(2s+1)+3\cdot2^{r+3}t+6=6(m+2^{r+2}t+1)+2^{r+1}(2s+1).
\end{equation}
Now (5.3) becomes,
\begin{equation}6(m+2^{r+2}t+1)+2^{r+1}(2s+1)\equiv 0\quad(\mbox{mod}\;2^{r+3}(2s+1));\end{equation}
equivalently,
\begin{equation}3(m+2^{r+2}t+1)+2^r(2s+1)\equiv 0\quad(\mbox{mod}\;2^{r+2}(2s+1)).\end{equation}

Note
\begin{equation}3(m+2^{r+2}t+1)+2^r(2s+1)\equiv 2^{r+2}(2s+1))\quad(\mbox{mod}\;2^{r+2}(2s+1));\end{equation}
equivalently,
\begin{equation}3(m+2^{r+2}t+1)\equiv 2^r(2s+1)3\quad(\mbox{mod}\;2^{r+2}(2s+1)).\end{equation}
Hence
\begin{equation}m+2^{r+1}t+1\equiv 2^r(2s+1)\quad(\mbox{mod}\;2^{r+1}(2s+1)),\end{equation}
that is,
\begin{equation}m=2^r(2s+1)+2^{r+2}(2s+1)c-2^{r+2}t-1\qquad\mbox{for some}\;\;c\in\mbb N.\end{equation}
Moreover,
\begin{equation}\Im_2=4\ell+3=3\cdot2^{r+5}t+23,\quad 6m+j+\ell+1=2^{r+3}(2s+1)(3c+1).\end{equation}

According to the second expression in (2.28) and (2.30),
\begin{equation}(3\cdot2^{r+5}t+23)|(2^{r+3}(2s+1)(3c+1))\lra (3\cdot2^{r+5}t+23)|(3c+1).\end{equation}
Thus
\begin{equation}c=2^{r+6}t+15+d(3\cdot2^{r+5}t+23)\qquad\mbox{for some}\;\;d\in\mbb N.\end{equation}
Now
\begin{equation}m=2^r(2s+1)[4(2^{r+6}t+15+d(3\cdot2^{r+5}t+23))+1]-2^{r+2}t-1.\end{equation}
and
\begin{equation} 6m+j+\ell+1=(2d+2)2^{r+3}(2s+1)(3\cdot2^{r+5}t+23)).\end{equation}
Therefore, (2.11) becomes
\begin{eqnarray}&&\frac{4}{24m+1}=\frac{1}{(3d+2)2^{r+3}(2s+1)(3\cdot2^{r+5}t+23)}\nonumber\\ & &+\frac{2^{r+3}(2s+1)+(3\cdot2^{r+5}t+23)}
{(3d+2)2^{r+3}(2s+1)(3\cdot2^{r+5}t+23)(24m+1)}\nonumber\\&=&
\frac{1}{(3d+2)2^{r+3}(2s+1)(3\cdot2^{r+5}t+23)}\nonumber\\ & &+\frac{1}{(3d+2)(3\cdot2^{r+5}t+23)(24m+1)}
\nonumber\\ & &+\frac{1}{(3d+2)2^{r+3}(2s+1)(24m+1)}.\end{eqnarray}
\pse

{\bf Theorem 5.2}\quad {\it  When $m$ is in (5.25), the Erd\"{o}s-Straus equation (5.27) holds.}\psp

\subsection{$j=2^{2r+2}(6s+1)$ and $\ell=3\cdot2^{2r+4}t+1$}

In this case,
\begin{equation}6m+j+\ell+1=6m+2^{2r+2}(6s+1)+3\cdot2^{2r+4}t+2=6[m+2^{2r+2}(s+2t)+(2^{2r+1}+1)/3].
\end{equation}
Now (5.3) becomes,
\begin{equation}6[m+2^{2r+2}(s+2t)+(2^{2r+1}+1)/3]\equiv 0\quad(\mbox{mod}\;2^{2r+4}(6s+1));\end{equation}
equivalently,
\begin{equation}3[m+2^{2r+2}(s+2t)+(2^{2r+1}+1)/3]\equiv 0\quad(\mbox{mod}\;2^{2r+3}(6s+1)).\end{equation}
So
\begin{equation}m+2^{2r+2}(s+2t)+(2^{2r+1}+1)/3\equiv 0\quad(\mbox{mod}\;2^{2r+3}(6s+1)).\end{equation}
Thus
\begin{equation}m=2^{2r+3}(6s+1)c-2^{2r+2}(s+2t)-(2^{2r+1}+1)/3\qquad\mbox{for some}\;\;0<c\in\mbb N.\end{equation}
Note
\begin{equation}\Im_2=4\ell+3=3\cdot2^{2r+6}t+7,\quad 6m+j+\ell+1=2^{2r+4}3(6s+1)c.\end{equation}
According to the second expression in (2.28) and (2.30),,
\begin{equation}(3\cdot2^{2r+6}t+7)|(3\cdot2^{2r+4}(6s+1)c)\end{equation}
\begin{equation}\lra c=(d+1)(3\cdot2^{r+6}t+7)\qquad\mbox{for some}\;\;d\in\mbb N.\end{equation}
Hence
\begin{equation}m=(d+1)2^{2r+3}(6s+1)(3\cdot2^{2r+6}t+7)-2^{2r+2}(s+2t)-(2^{2r+1}+1)/3\end{equation}
and
\begin{equation} 6m+j+\ell+1=3(d+1)2^{2r+4}(6s+1)(3\cdot2^{2r+6}t+7)d.\end{equation}
Therefore, (2.11) becomes
\begin{eqnarray}&&\frac{4}{24m+1}=\frac{1}{3(d+1)2^{2r+4}(6s+1)(3\cdot2^{2r+6}t+7)}\nonumber\\ & &+\frac{2^{2r+4}(6s+1)+(3\cdot2^{2r+6}t+7)}
{(d+1)2^{2r+4}(6s+1)(3\cdot2^{2r+6}t+7)(24m+1)}\nonumber\\&=&
\frac{1}{3(d+1)2^{2r+4}(6s+1)(3\cdot2^{2r+6}t+7)}\nonumber\\&&+\frac{1}{3(d+1)(3\cdot2^{2r+6}t+7)(24m+1)}
\nonumber\\ & &+\frac{1}{3(d+1)2^{2r+4}(6s+1)(24m+1)}.\end{eqnarray}

{\bf Theorem 5.3}\quad {\it  When $m$ is in (5.36), the Erd\"{o}s-Straus equation (5.38) holds}.\psp

\subsection{$j=2^{2r+2}(6s+5)$ and $\ell=3\cdot2^{2r+4}t+1$}

In this case,
\begin{equation}6m+j+\ell+1=6m+2^{2r+2}(6s+5
)+3\cdot2^{2r+4}t+2=6(m+2^{2r+3}t))+2^{2r+2}(6s+5)+2.
\end{equation}
Now (5.3) becomes,
\begin{equation}6(m+2^{2r+3}t))+2^{2r+2}(6s+5)+2\equiv 0\quad(\mbox{mod}\;2^{2r+4}(6s+5));\end{equation}
equivalently,
\begin{equation}3(m+2^{2r+3}t))+2^{2r+1}(6s+5)+1\equiv 0\quad(\mbox{mod}\;2^{2r+3}(6s+5)).\end{equation}

Note
\begin{equation}3(m+2^{2r+3}t))+2^{2r+1}(6s+5)+1\equiv 2^{2r+4}(6s+5)\quad(\mbox{mod}\;2^{2r+3}(6s+5)).\end{equation}
So
\begin{equation}3(m+2^{2r+3}t))\equiv 2^{2r+1}7(6s+5)-1\quad(\mbox{mod}\;2^{2r+3}(6s+5)).\end{equation}
\begin{equation}3(m+2^{2r+3}t)\equiv 3[2^{2r+2}(7s+6)-(2^{2r+1}+1)/3]\quad(\mbox{mod}\;2^{2r+3}(6s+5)).\end{equation}
Hence
\begin{equation}m+2^{2r+3}t\equiv 2^{2r+2}(7s+6)-(2^{2r+1}+1)/3\quad(\mbox{mod}\;2^{2r+3}(6s+5)).\end{equation}
Thus
\begin{equation}m=2^{2r+3}(6s+5)c+2^{2r+2}(7s-2t+6)-(2^{2r+1}+1)/3\qquad\mbox{for some}\;\;c\in\mbb N.\end{equation}

Observe
\begin{equation}\Im_2=4\ell+3=3\cdot2^{2r+6}t+7,\quad 6m+j+\ell+1=(3c+2)2^{2r+4}(6s+5).\end{equation}
According to the second expression in (2.28) and (2.30),
\begin{equation}(3\cdot2^{2r+6}t+7)|((3c+2)2^{2r+4}(6s+5))\lra (3\cdot2^{2r+6}t+7)|(3c+2)  \end{equation}
\begin{equation}\lra c=2^{2r+7}t+4+d(3\cdot2^{2r+6}t+7)\qquad\mbox{for some}\;\;d\in\mbb N.\end{equation}
Hence
\begin{equation}m=2^{2r+3}(6s+5)(2^{2r+7}t+4+d(3\cdot2^{2r+6}t+7))+2^{2r+2}(7s-2t+6)-(2^{2r+1}+1)/3\end{equation}
and
\begin{equation} 6m+j+\ell+1=(3d+2)2^{2r+4}(6s+5)(3\cdot2^{2r+6}t+7).\end{equation}
Therefore, (2.11) becomes
\begin{eqnarray}&&\frac{4}{24m+1}=\frac{1}{(3d+2)2^{2r+4}(6s+5)(3\cdot2^{2r+6}t+7)}\nonumber\\ & &+\frac{2^{2r+4}(6s+5)+(3\cdot2^{2r+6}t+7)}{(3d+2)2^{2r+4}(6s+5)(3\cdot2^{2r+6}t+7)(24m+1)}\nonumber\\&=&
\frac{1}{(3d+2)2^{2r+4}(6s+5)(3\cdot2^{2r+6}3t+7)}\nonumber\\ & &+\frac{1}{(3d+2)(3\cdot2^{2r+6}t+7)(24m+1)}
\nonumber\\ & &+\frac{1}{(3d+2)2^{2r+4}(6s+5)(24m+1)}.\end{eqnarray}\pse

{\bf Theorem 5.4}\quad {\it  When $m$ is in (5.50), the Erd\"{o}s-Straus equation (5.52) holds}.\psp

\subsection{$j=2^{2r+2}(6s+5)$ and $\ell=3\cdot2^{2r+4}t+3$}

In this case,
\begin{equation}6m+j+\ell+1=6m+2^{2r+2}(6s+5)+3\cdot2^{2r+4}t+4=6[m+2^{2r+2}(s+2t)+(5\cdot 2^{2r+1}+2)/3].
\end{equation}
Now (5.3) becomes,
\begin{equation}6[m+2^{2r+2}(s+2t)+(\cdot 2^{2r+1}+2)/3]\equiv 0\quad(\mbox{mod}\;2^{2r+4}(6s+5));\end{equation}
equivalently,
\begin{equation}3[m+2^{2r+2}(s+2t)+(5\cdot2^{2r+1}+2)/3]\equiv 0\quad(\mbox{mod}\;2^{2r+3}(6s+5)).\end{equation}
So
\begin{equation}m+2^{2r+2}(s+2t)+(5\cdot2^{2r+1}+2)/3\equiv 0\quad(\mbox{mod}\;2^{2r+3}(6s+5)).\end{equation}
Thus
\begin{equation}m=2^{2r+3}(6s+5)c-2^{2r+2}(s+2t)-(5\cdot2^{2r+1}+2)/3\qquad\mbox{for some}\;\;0<c\in\mbb N.\end{equation}
Note
\begin{equation}\Im_2=4\ell+3=3\cdot2^{2r+6}t+15,\quad 6m+j+\ell+1=3\cdot2^{2r+4}(6s+5)c.\end{equation}

According to the second expression in (2.28) and (2.30),
\begin{equation}(2^{2r+6}3t+15)|(3\cdot2^{2r+4}(6s+1)c)\end{equation}
\begin{equation}\lra c=(d+1)(2^{r+6}t+5)\qquad\mbox{for some}\;\;d\in\mbb N.\end{equation}
Hence
\begin{equation}m=(d+1)2^{2r+3}(6s+5)(2^{2r+6}t+5)d-2^{2r+2}(s+2t)-(5\cdot2^{2r+1}+2)/3\end{equation}
and
\begin{equation} 6m+j+\ell+1=(d+1)2^{2r+4}(6s+5)(3\cdot2^{2r+4}t+15)d.\end{equation}
Therefore, (2.11) becomes
\begin{eqnarray}&&\frac{4}{24m+1}=\frac{1}{(d+1)2^{2r+4}(6s+5)(3\cdot2^{2r+6}t+15)}\nonumber\\ & &+\frac{2^{2r+4}(6s+5)+(3\cdot2^{2r+6}t+15)}
{(d+1)2^{2r+4}(6s+5)(3\cdot2^{2r+6}t+15)(24m+1)}\nonumber\\&=&
\frac{1}{(d+1)2^{2r+4}(6s+1)(3\cdot2^{2r+6}t+15)}\nonumber\\ & &+\frac{1}{(d+1)(3\cdot2^{2r+6}t+15)(24m+1)}
\nonumber\\ & &+\frac{1}{(d+1)2^{2r+4}(6s+5)d(24m+1)}.\end{eqnarray}\pse

{\bf Theorem 5.5
}\quad {\it  When $m$ is in (5.61), the Erd\"{o}s-Straus equation (5.63) holds}.\psp

\subsection{$j=2^{2r+1}(6s+1)$ and $\ell=3\cdot2^{2r+3}t+1$}

In this case,
\begin{equation}6m+j+\ell+1=6m+2^{2r+1}(6s+1)+3\cdot2^{2r+3}t+2=6(m+2^{2r+2}t)+2^{2r+1}(6s+1)+2.
\end{equation}
Now (5.3) becomes,
\begin{equation}6(m+2^{2r+2}t)+2^{2r+1}(6s+1)+2\equiv 0\quad(\mbox{mod}\;2^{2r+3}(6s+1));\end{equation}
equivalently,
\begin{equation}3(m+2^{2r+2}t)+2^{2r}(6s+1)+1\equiv 0\quad(\mbox{mod}\;2^{2r+2}(6s+1)).\end{equation}

Note
\begin{equation}3(m+2^{2r+2}t)+2^{2r}(6s+1)+1\equiv 2^{2r+3}(6s+1)\quad(\mbox{mod}\;2^{2r+2}(6s+1))\end{equation}
\begin{equation}\lra 3(m+2^{2r+2}t)\equiv 2^{2r}7(6s+1)-1\quad(\mbox{mod}\;2^{2r+2}(6s+1))\end{equation}
\begin{equation}\lra 3(m+2^{2r+2}t)\equiv 3[2^{2r+1}(7s+1)+(2^{2r}-1)/3]\quad(\mbox{mod}\;2^{2r+2}(6s+1))\end{equation}
\begin{equation}\lra m+2^{2r+2}t\equiv 2^{2r+1}(7s+1)+(2^{2r}-1)/3.\quad(\mbox{mod}\;2^{2r+2}(6s+1))\end{equation}
So
\begin{equation}m=c(2^{2r+2}(6s+1))+2^{2r+1}(7s-2t+1)+(2^{2r}-1)/3\qquad\mbox{for some}\;\;c\in\mbb N.\end{equation}

Now
\begin{equation}\Im_2=4\ell+3=3\cdot2^{2r+5}t+7,\quad 6m+j+\ell+1=3(3c+2)2^{2r+3}(6s+1).\end{equation}
According to the second expression in (2.28) and (2.30),
\begin{equation}(3\cdot2^{2r+5}t+7)|(3(3c+2)2^{2r+3}(6s+1))\lra(3\cdot2^{2r+5}t+7)|(3c+2) \end{equation}
\begin{equation}\lra c=2^{2r+6}t+4+d(3\cdot2^{2r+5}t+7)\qquad\mbox{for some}\;\;d\in\mbb N.\end{equation}
Hence
\begin{equation}m=2^{2r+2}(6s+1)(2^{2r+6}t+4+d(3\cdot2^{2r+5}t+7))+2^{2r+1}(7s-2t+1)+(2^{2r}-1)/3\end{equation}
and
\begin{equation} 6m+j+\ell+1=3(3d+2)2^{2r+3}(6s+1)(3\cdot2^{2r+5}t+7).\end{equation}
Therefore, (2.11) becomes
\begin{eqnarray}&&\frac{4}{24m+1}=\frac{1}{3(3d+2)2^{2r+3}(6s+1)(3\cdot2^{2r+5}t+7)}\nonumber\\ & &+\frac{2^{2r+3}(6s+1)+(3\cdot2^{r+5}t+7)}
{3(3d+2)2^{2r+3}(6s+1)(3\cdot2^{2r+5}t+7)(24m+1)}\nonumber\\&=&
\frac{1}{3(3d+2)2^{2r+3}(6s+1)(3\cdot2^{r+5}t+7)}\nonumber\\ & &+\frac{1}{3(3d+2)(3\cdot2^{2r+5}t+7)(24m+1)}
\nonumber\\ & &+\frac{1}{3(3d+2)2^{2r+3}(6s+1)(24m+1)}.\end{eqnarray}\pse

{\bf Theorem 5.6
}\quad {\it  When $m$ is in (5.75), the Erd\"{o}s-Straus equation (5.77) holds}.\psp

\subsection{$j=2^{2r+1}(6s+1)$ and $\ell=3\cdot2^{2r+3}t+3$}

In this case,
\begin{equation}6m+j+\ell+1=6m+2^{2r+1}(6s+1)+3\cdot2^{2r+3}t+4=6[m+2^{2r+1}(s+2t)+(2^{2r}+2)/3].
\end{equation}
Now (5.3) becomes,
\begin{equation}6[m+2^{2r+1}(s+2t)+(2^{2r}+2)/3]\equiv 0\quad(\mbox{mod}\;2^{2r+3}(6s+1));\end{equation}
equivalently,
\begin{equation}3[m+2^{2r+1}(s+2t)+(2^{2r}+2)/3]\equiv 0\quad(\mbox{mod}\;2^{2r+2}(6s+1)).\end{equation}
So
\begin{equation}m+2^{2r+1}(s+2t)+(2^{2r}+2)/3\equiv 0\quad(\mbox{mod}\;2^{2r+2}(6s+1)).\end{equation}
Thus
\begin{equation}m=2^{2r+2}(6s+1)c-2^{2r+1}(s+2t)-(2^{2r}+2)/3.\qquad\mbox{for some}\;\;c\in\mbb N.\end{equation}
Note
\begin{equation}\Im_2=4\ell+3=3\cdot2^{2r+5}t+15,\quad 6m+j+\ell+1=3\cdot2^{2r+3}(6s+1)c.\end{equation}
According to the second expression in (2.28) and (2.30),
\begin{equation}(3\cdot2^{2r+5}t+15)|(3\cdot2^{2r+3}(6s+1)c)\end{equation}
\begin{equation}\lra c=(d+1)(2^{r+5}t+5)\qquad\mbox{for some}\;\;d\in\mbb N.\end{equation}
Hence
\begin{equation}m=(d+1)2^{2r+2}(6s+1)(2^{2r+5}t+5)-2^{2r+1}(s+2t)-(2^{2r}+2)/3\end{equation}
and
\begin{equation} 6m+j+\ell+1=(d+1)2^{2r+3}(6s+1)(3\cdot2^{2r+5}t+15).\end{equation}
Therefore, (2.11) becomes
\begin{eqnarray}&&\frac{4}{24m+1}=\frac{1}{(d+1)2^{2r+3}(6s+1)(3\cdot2^{2r+5}t+15)}\nonumber\\ & &+\frac{2^{2r+3}(6s+1)+(3\cdot2^{r+5}t+15)}
{(d+1)2^{2r+3}(6s+1)(3\cdot2^{2r+5}t+15)(24m+1)}\nonumber\\&=&
\frac{1}{(d+1)2^{2r+3}(6s+1)(3\cdot2^{r+5}t+15)}\nonumber\\ & &+\frac{1}{(d+1)(3\cdot2^{r+5}t+15)(24m+1)}
\nonumber\\ & &+\frac{1}{(d+1)2^{2r+3}(6s+1)d(24m+1)}.\end{eqnarray}\pse

{\bf Theorem 5.7}\quad {\it  When $m$ is in (5.86), the Erd\"{o}s-Straus equation (5.88) holds}.

\subsection{$j=2^{2r+1}(6s+5)$ and $\ell=3\cdot2^{2r+3}t+1$}

In this case,
\begin{equation}6m+j+\ell+1=6m+2^{2r+1}(6s+5)+3\cdot2^{2r+3}t+2=6[m+2^{2r+1}(s+2t)+(5\cdot2^{2r}+1)/3].
\end{equation}
Now (5.3) becomes,
\begin{equation}6[m+2^{2r+1}(s+2t)+(5\cdot2^{2r}+1)/3]\equiv 0\quad(\mbox{mod}\;2^{2r+3}(6s+5));\end{equation}
equivalently,
\begin{equation}3[m+2^{2r+1}(s+2t)+(5\cdot2^{2r}+1)/3]\equiv 0\quad(\mbox{mod}\;2^{2r+2}(6s+5)).\end{equation}
So
\begin{equation}m+2^{2r}(s+2t)+(5\cdot2^{2r}+1)/3\equiv 0\quad(\mbox{mod}\;2^{2r+2}(6s+5)).\end{equation}
Thus
\begin{equation}m=2^{2r+2}(6s+5)c-2^{2r+1}(s+2t)-(5\cdot2^{2r}+1)/3\qquad\mbox{for some}\;\;c\in\mbb N.\end{equation}
Note
\begin{equation}\Im_2=4\ell+3=3\cdot2^{2r+5}t+7,\quad 6m+j+\ell+1=3\cdot2^{2r+3}(6s+5)c.\end{equation}
According to the second expression in (2.28) and (2.30),
\begin{equation}(3\cdot2^{2r+5}t+7)|(3\cdot2^{2r+3}(6s+1)c)\end{equation}
\begin{equation}\lra c=(d+1)(3\cdot2^{2r+5}t+7)\qquad\mbox{for some}\;\;d\in\mbb N.\end{equation}
Hence
\begin{equation}m=(d+1)2^{2r+2}(6s+5)(3\cdot2^{2r+5}t+7)-2^{2r+1}(s+2t)-(5\cdot2^{2r}+1)/3\end{equation}
and
\begin{equation} 6m+j+\ell+1=3(d+1)2^{2r+3}(6s+5)(3\cdot2^{2r+5}t+7).\end{equation}
Therefore, (2.11) becomes
\begin{eqnarray}&&\frac{4}{24m+1}=\frac{1}{3(d+1)2^{2r+3}(6s+5)(3\cdot2^{2r+5}t+7)}\nonumber\\ & &+\frac{2^{2r+3}(6s+5)+(3\cdot2^{r+5}t+7)}
{3(d+1)2^{2r+3}(6s+5)(3\cdot2^{2r+5}t+7)(24m+1)}\nonumber\\&=&
\frac{1}{3(d+1)2^{2r+3}(6s+5)(3\cdot2^{r+5}t+7)}+\frac{1}{3(d+1) (3\cdot2^{2r+5}t+7)(24m+1)}
\nonumber\\ & &+\frac{1}{3(d+1)2^{2r+3}(6s+5)(24m+1)}.\end{eqnarray}\pse

{\bf Theorem 5.8
}\quad {\it  When $m$ is in (5.97), the Erd\"{o}s-Straus equation (5.99) holds}.

\subsection{Section Summary}

In the above sections, $m=g(x,d,y)$ and $\Im_2=h(y)$ are integral polynomials, where  $y=2^{r+1}t$, $y=2^{2+1}t$ or $y=2^{2r+2}t$.
If we change $y$ to any positive integer, they still give a solution of the Erd\"{o}s-Straus equation. In summary, we have the following table
of solutions:

\psp
\begin{center}{\bf \large Table 5.6}\end{center}\pse
{\small \begin{center}\begin{tabular}{|c|c|c|}\hline $m$&$\Im_1$&$\Im_2$\\\hline $(d+1)2^{r+2}(2x+1)(24y+23)-y-2^{r}(2x+1)-1$&
$2^{r+3}(6x+3)$&$24y+23$\\\hline $2^r(2x+1)[4(16y+15+d(24y+23))+1]-y-1$&$2^{r+3}(2x+1)$&$24y+23$ \\\hline
$(d+1)2^{2r+3}(6x+1)(24y+7)-2y-2^{2r+2}x-(2^{2r+1}+1)/3$&$2^{2r+4}(6x+1)$&$24y+7$\\\hline
$2^{2r+2}[2(6x+5)(16y+4+d(24y+7))+7x+6]-y-(2^{2r+1}+1)/3$&$2^{2r+4}(6x+5)$&$24y+7$\\\hline
$(d+1)2^{2r+3}(6x+5)(8y+5)-y-2^{2r}x-(5\cdot2^{2r+1}+2)/3$&$2^{2r+4}(6x+5)$&$24y+15$\\\hline
$2^{2r+1}[2(6x+1)(16y+4+d(24y+7))+7x+1]-y+(2^{2r}-1)/3$&$2^{2r+3}(6x+1)$&$24y+7$\\\hline
$(d+1)2^{2r+2}(6x+1)(8y+5)-y-2^{2r+1}x-(2^{2r}+2)/3$&$2^{2r+3}(6x+1)$&$24y+15$\\\hline
$(d+1)2^{2r+2}(6x+5)(24y+7)-y-2^{2r+1}x-(3\cdot2^{2r}+1)/3$&$2^{2r+3}(6x+5)$&$24y+7$\\\hline
\end{tabular}\end{center}}

\psp

\psp

\psp

E-Mail:  xiaoping@math.ac.cn
\pse

 \end{document}